\documentclass[12pt]{article}

\usepackage{latexsym}
\usepackage{amsfonts}
\usepackage{amssymb}

\newcommand{\so}[1]{\mbox{${\mathfrak s \mathfrak o}($}}

\newcommand{\g}{\mbox{${\mathfrak g}$}}
\newcommand{\h}{\mbox{${\mathfrak h}$}}

\newcommand{\p}{\mbox{${\mathfrak p}$}}

\newcommand{\cl}{\mbox{${\mathbb C}$}}

\newcommand{\rl}{\mbox{${\mathbb R}$}}

\makeatletter
\def\numberwithin#1#2{\@ifundefined{c@#1}{\@nocnterrr}{%
  \@ifundefined{c@#2}{\@nocnterr}{%
  \@addtoreset{#1}{#2}%
  \toks@\expandafter\expandafter\expandafter{\csname the#1\endcsname}%
  \expandafter\xdef\csname the#1\endcsname
    {\expandafter\noexpand\csname the#2\endcsname
     .\the\toks@}}}}
\makeatother
\numberwithin{equation}{section}

\newtheorem{theorem}[equation]{Theorem}
\newtheorem{lemma}[equation]{Lemma}

\newtheorem{ex}[equation]{Example}

\newtheorem{rem}[equation]{Remark}
\newenvironment{remark}{\begin{rem} \em}{\end{rem}}

\begin{document}

\author{Andrew S. Dancer and  Ian A. B. Strachan}
\title{ \bf Einstein metrics on tangent bundles of spheres.}
\date{\small \today}
\maketitle

\setcounter{section}{0}
\setcounter{page}{1}
\setcounter{equation}{0}

\begin{abstract}
We give an elementary treatment of the existence of complete K\"ahler-Einstein
metrics with nonpositive Einstein constant and underlying manifold
diffeomorphic to the tangent bundle of the $(n+1)$-sphere.
\end{abstract}

\noindent{ {\it Mathematics Subject Classification (2000) 53C}

\medskip
\section{Introduction.}
Over the last few years there has been considerable interest in
a family of Ricci-flat K\"ahler metrics discovered by
Stenzel \cite{St} with underlying manifold
diffeomorphic to the tangent bundle of the $(n+1)$-sphere. (The $n=2$
case was found earlier by Candelas and de la Ossa \cite{CD} : the
$n=1$ case is the Eguchi-Hanson metric \cite{EH}).
The  complex structure is that of the quadric in $\cl^{n+2}$.

One source of interest in physics is that there is a ``conifold
transition''\cite{CD}--both the Stenzel metrics and another family of
Ricci-flat K\"ahler metrics due to B\'erard-Bergery \cite{BB} have a
cone as a common degenerate limit. Also, Cvetic, Gibbons, L\"u and
Pope have recently studied the harmonic forms on these metrics and
found an explicit formula for the Stenzel metrics in terms of
hypergeometric functions \cite{CGLP}.

In an unpublished 1993 preprint \cite{DS2} we gave a rather different
description of these metrics, and showed that they also have analogues
with negative Einstein constant. The underlying manifold of the latter
metrics is still diffeomorphic to $TS^{n+1}$, though the complex
structure is now that of a tube in the quadric rather than the whole
quadric.

As several people have expressed interest in this work, we hope that
it may be worthwhile to publish it, with some amplification,
 in the current paper.

We remark that the existence of a canonical complex structure on (a
tube in) the tangent bundle of a Riemannian manifold $(M, g)$ is a
general fact due to Guillemin-Stenzel \cite{GS} and Lempert-Sz\"oke
\cite{LSz}, \cite{Sz1}.  Now the existence theorem of Mok-Yau \cite{MY}
shows the existence of complete K\"ahler-Einstein metrics on tubes in
$TM$ of radius smaller than the maximal radius on which the complex
structure is defined \cite{Sz2}.

In our case (i.e. when $M$ is the round sphere) we are able to give a
much more elementary and concrete treatment,  exploiting the
fact that the metrics are of cohomogeneity one with respect to an
action of $SO(n+2)$.  Our techniques are similar to those of our paper
\cite{DS1}, where we classified K$\rm \ddot{a}$hler-Einstein metrics
in real dimension four which are of Bianchi IX type.

\bigskip
\section{Cohomogeneity one K\"ahler-Einstein metrics. }

\bigskip
Let us consider a Riemannian manifold $(M,g)$, which admits an isometric action
of $SO(n+2)$ with principal orbit $SO(n+2)/SO(n)$. We further suppose
that this action is of cohomogeneity one, that is, each principal orbit
has real codimension one in $M$. It follows that the (real) dimension
of $M$ is $2n+2$.

The union of the principal orbits will form an open dense set in $M$.
On this set we can write the metric as
\[
g=dt^{2} + g_{t}
\]
where $t$ is the arclength parameter along a geodesic orthogonal
to the group orbits, and $g_{t}$ is a homogeneous metric on the orbits.

Now, for any homogeneous space $G/H$ with $H$ compact we may identify
the tangent space at a point with an Ad $H$ -invariant complement $\p$
for $\h$ in $\g$. With this identification, $G$-invariant metrics on
$G/H$ correspond precisely to Ad $H$-invariant inner products on $\p$.

In our case we embed $\so (n)$ in $\so(n+2)$ so that if $X_{i,j}$
 denotes the matrix with $ij$th entry $=1$, $ji$th entry= $-1$ and all
 other entries zero, then $\{X_{i,j}  : \; j=n+1,n+2 ,1 \leq i <j \}$ spans
 a complement $\p$ for $\so(n)$ in $\so(n+2)$.  Now the decomposition
 of $\p$ into irreducible components under the adjoint action of
 $SO(n)$ is
\[
\p= \p_{1} \oplus \p_{2} \oplus \p_{3}
\]
where $\p_{1},\p_{2}$ are standard $n$-dimensional representations
 spanned by $\{X_{i,n+1} : 1 \leq i \leq n \}$ and $\{X_{i,n+2} : 1
 \leq i \leq n \}$ respectively, and $\p_{3}$ is the trivial
 representation spanned by $X_{n+1,n+2}$.

We shall consider metrics of the form
\[
g_t =  a(t)^{2} \;B \mid_{\p_{1}} \oplus \; b(t)^{2} \;B
\mid_{\p_{2}} \oplus \; c(t)^{2} \;B \mid_{\p_{3}}
\]
where $B$ is the Killing form on $\so(n+2)$ defined by
$B(X,Y)= -\frac{1}{2} {\rm Tr} \;XY$.
Note that this is not the most general possibility for $g_t$. As
$\p_1$ and $\p_2$ are isomorphic as $SO(n)$-representations it is a
definite restriction to assume that $\p_1$ and $\p_2$ are orthogonal
for all $t$.

If $a^{2}=b^{2}$ our metric $g_t$ is obtained using a Riemannian
submersion with circle fibre over the flag manifold $SO(n+2)/SO(n)
\times SO(2)$ equipped with a K$\rm \ddot{a}$hler-Einstein metric.  If
$a^{2}$ and $b^{2}$ are not equal, however, $g_{t}$ does not arise in
this way.

It is convenient to introduce at this stage a new transverse
coordinate $\zeta$ defined by
\[
dt = a^{n}b^{n}c \; d \zeta.
\]
We shall denote differentiation with respect to $\zeta$ by a prime.

The formulae of B$\rm \acute{e}$rard-Bergery \cite{BB}
 enable us to write down the Einstein equations
${\mbox Ric}_g = \Lambda g$ for the cohomogeneity one metric $g$. They are:
\begin{eqnarray}
-\frac{1}{a^{2n-2}b^{2n}c^{2}} \left( \frac{a^{\prime}}{a} \right)^{\prime}
+ \frac{a^{4} -(b^{2}-c^{2})^{2}}{2b^{2}c^{2}} +n-1 &=& \Lambda a^{2},
\label{BBE1}\\
-\frac{1}{a^{2n}b^{2n-2}c^{2}} \left( \frac{b^{\prime}}{b} \right)^{\prime}
+ \frac{b^{4} -(c^{2}-a^{2})^{2}}{2a^{2}c^{2}} + n-1 &=& \Lambda
b^{2}, \label{BBE2}\\
-\frac{1}{a^{2n}b^{2n}} \left( \frac{c^{\prime}}{c} \right)^{\prime}
+ \frac{n(c^{4}-(a^{2}-b^{2})^{2})}{2a^{2}b^{2}} &=& \Lambda c^{2},
\label{BBE3}
\end{eqnarray}

\begin{equation}
-\frac{1}{a^{n}b^{n}c} \left(\frac{1}{a^{n}b^{n}c} \left(
 n \frac{a^{\prime}}{a} + n \frac{b^{\prime}}{b}+
 \frac{c^{\prime}}{c} \right) \right)^{\prime}
-\frac{1}{a^{2n}b^{2n}c^{2}} \left(
 n \left(\frac{a^{\prime}}{a} \right)^{2}
 + n \left(\frac{b^{\prime}}{b} \right)^{2}
+ \left(\frac{c^{\prime}}{c} \right)^{2}  \right) = \Lambda, \label{BBE4}
\end{equation}
where $\Lambda$ is the Einstein constant.

\begin{remark}
If $n=1$, then (\ref{BBE1})-(\ref{BBE4}) are a presentation of the
 Einstein equations for a Bianchi IX metric.
\end{remark}

It can be verified by direct calculation that
(\ref{BBE1})-(\ref{BBE4}) hold in particular if $a,b,c$ satisfy the
following three first-order equations
\begin{eqnarray}
a^{\prime} &=& \frac{1}{2} a^{n}b^{n-1}(b^{2}+c^{2}-a^{2}), \label{BBKE1} \\
b^{\prime} &=& \frac{1}{2}a^{n-1}b^{n}(c^{2}+a^{2}-b^{2}), \label{BBKE2}\\
c^{\prime} &=& \frac{1}{2}na^{n-1}b^{n-1}c \left( a^{2}+b^{2}-c^{2} -
\frac{2 \Lambda}{n}a^{2}b^{2} \right). \label{BBKE3}
\end{eqnarray}
It is this reduction of the Einstein equations which we shall study.
The metrics arising from solutions of these equations are precisely
those Einstein metrics with respect to which the $SO(n+2)$-invariant
almost complex structure $J$ defined by
\begin{equation} \label{cplex}
J : X_{i,n+1} \mapsto -\frac{a}{b} X_{i, n+2}, \;\;\; : \;\;\;
\frac{\partial}{\partial t} \mapsto -\frac{1}{c} X_{n+1,n+2}
\end{equation}
 is K$\rm \ddot{a}$hler.

\bigskip
\section{Complete metrics.}

\bigskip
We shall now demonstrate the existence of solutions to
(\ref{BBKE1})-(\ref{BBKE3}) which give rise to complete metrics.  It
will be convenient to introduce a new variable $u$ defined by
\[
du =(ab)^{n-1} d \zeta
\]
so $dt = abc \; du$, and our equations become
\begin{eqnarray}
a_{u} &=& \frac{1}{2}a(b^{2}+c^{2}-a^{2}), \label{KE1} \\
b_{u} &=& \frac{1}{2}b(c^{2}+a^{2}-b^{2}), \label{KE2} \\
c_{u} &=& \frac{1}{2}nc \left(a^{2}+b^{2}-c^{2}
-\frac{2 \Lambda}{n}a^{2}b^{2} \right), \label{KE3}.
\end{eqnarray}

\begin{remark}
If $n=1$ these are the Bianchi IX K\"ahler-Einstein
equations studied by the authors in \cite{DS1}.
\end{remark}

The critical points of the system (\ref{KE1})-(\ref{KE3})
are precisely the points
$(0,K,\pm K),(K,0, \pm K)$ and $(K, \pm K,0)\,,$ where $K \in \rl$.

\medskip
{\bf Assumption}

{}From now on we shall take $\Lambda$ to be less than or equal to zero.

\bigskip
Subject to this assumption, the linearisation of (\ref{KE1})-(\ref{KE3}) about
 each critical point (except the origin)
has one negative, one positive and one zero eigenvalue.
So each of these critical points has an unstable curve.

Let us consider an unstable curve of $(0,K,K)$ where $K$ is nonzero.
This will be a solution to (\ref{KE1})-(\ref{KE3}) defined on a
maximal interval $(-\infty, \eta)$ for some $\eta$ ($\eta$ may be
$\infty$).

It follows from (\ref{KE1})-(\ref{KE3}) that if any one of $a,b$ or
$c$ is zero at a point in $(-\infty, \eta)$, then it is identically
zero. It is clear that on the unstable curve none of $a,b$ or $c$ is
identically zero, so none of them vanishes anywhere on $(-\infty,
\eta)$. The metric and the equations
 (\ref{KE1})-(\ref{KE3}) are invariant under
changes of sign of $a,b$ or $c$ so from now on we can assume that
$a,b,c$ are all positive on $(-\infty, \eta)$; in particular we can
take $K$ to be positive.

 The metric is therefore
defined on $(-\infty, \eta)$ and to show that it is complete we need
to study its behaviour as $u$ tends to $-\infty$ and as $u$ tends to $\eta$.

Note that if $a$ equals $b$ at any point on the unstable curve
 then $a$ is identically equal to $b$,
giving a contradiction. It follows that $b$ is greater than $a$ on the
 unstable curve, and so (from (\ref{KE1})) $a$ is strictly increasing. It now
follows from (\ref{KE1}),(\ref{KE2}) that $\frac{b}{a}$
 is strictly decreasing, and
tends to some limit $L \geq 1$ as $u$ tends to $\eta$.

\begin{lemma}
On an unstable curve of $(0,K,K)$ we have the inequalities
\begin{equation}
c^{2} \leq a^{2}+b^{2}-\frac{2 \Lambda}{n} a^{2}b^{2} \label{lem1}
\end{equation}
\begin{equation}
b^{2} \leq a^{2}+c^{2}. \label{lem2}
\end{equation}
\end{lemma}

{\bf Proof}

\medskip
Suppose that $c^{2} > a^{2}+b^{2} -\frac{2 \Lambda}{n}a^{2}b^{2}$ at
$u_{0}$. The equations (\ref{KE2}),(\ref{KE3})
imply that $b$ is increasing and $c$ is
decreasing at $u_{0}$. Recall also that $a$ is increasing.
We deduce that $c^{2}$ is greater than
$a^{2}+b^{2}-\frac{2 \Lambda}{n}a^{2}b^{2}$ on $(-\infty, u_{0})$
and hence that $a,b$ are increasing and $c$ is decreasing on this interval.
It follows
 that $c$ is bounded away from $b$ on $(-\infty, u_{0})$.
This contradicts the fact that $(a,b,c)$ tends to $(0,K,K)$ as $u$
tends to $-\infty$, so we have established inequality (\ref{lem1}).
The proof of inequality (\ref{lem2}) is very similar. $\Box$

\begin{remark}
We deduce from (\ref{KE2}),(\ref{KE3}) that $b,c$ are
increasing on $(-\infty,\eta)$. We remarked earlier
that $a$ is strictly increasing on this interval.
\end{remark}

Let us now study the behaviour of the metric as $u$ tends to $-\infty$ and
as $u$ tends to $\eta$.

As $u$ tends to $-\infty$
\[
a \simeq M e^{K^{2}u}, \;\;\;
b \simeq K, \;\;\;
c \simeq K
\]
 for some constant $M$.

Choosing $R= Me^{K^{2}u}$ as a new coordinate,
find that the metric is asymptotically given by
\[
dR^{2} + R^{2} B \mid_{\p_{1}}
+ K^{2}( B \mid_{\p_{1}} \oplus B \mid_{\p_{3}}).
\]
as $R$ tends to zero.

Now, $B \mid_{\p_{1}}$ defines the standard $SO(n+1)$-invariant metric on
$S^{n}$, while $B \mid_{\p_{2}} \oplus B \mid_{\p_{3}}$ defines the standard
$SO(n+2)$-invariant metric on $S^{n+1}$.  Therefore we obtain a
nonsingular metric by adding in an $(n+1)$-sphere at $R=0$.  In terms
of the orbit type of the $SO(n+2)$ action, the isotropy group jumps at
$R=0$ from $SO(n)$ to $SO(n+1)$, so an $n$-dimensional sphere
collapses to a point and the orbit at $R=0$ is the $(n+1)$-sphere
$SO(n+2)/SO(n+1)$ rather than $SO(n+2)/SO(n)$. The orbit at $R=0$ is
called a {\em Bolt of the second kind} in the terminology of Gibbons,
Page and Pope \cite{GPP}. Note that $K$ determines and is determined
by  the volume of the Bolt.

To examine the behaviour of the metric as $u$ tends to $\eta$ we introduce
a new coordinate $\rho$, defined by
\[
\rho= 2 (ab)^{\frac{1}{2}}.
\]
This is an allowable change of variables as $ab$ is increasing.

The metric is now given by
\begin{equation}
W^{-1} d \rho^{2} + \frac{1}{4}\rho^{2}
\left( \frac{a}{b} \; B \mid_{\p_1} +  \frac{b}{a} \; B \mid_{\p_2}
+ W \; B \mid_{\p_3} \right) \label{metrho}
\end{equation}
 where $W=\frac{c^{2}}{ab}$.
It follows from equations (\ref{KE1})-(\ref{KE3}) that
\begin{equation}
\frac{dW}{d \rho} + \frac{2(n+1)W}{\rho} = \frac{2n}{\rho} \left( \frac{a}{b} +
 \frac{b}{a} \right) - \Lambda \rho. \label{eqrho}
\end{equation}
Recall that the inequalities (\ref{lem1}) and (\ref{lem2}) show that
$a,b,c$ are monotonic increasing on $(-\infty, \eta)$.  Suppose that
the limit $\lambda$ of $a$ as $u$ tends to $\eta$ is finite.  Since
$b/a$ is decreasing, and because of the estimate (\ref{lem1}), we see
that the limits $\mu,\nu$ of $b,c$ at $u=\eta$ are also finite.  If
$\eta=\infty$ then $(\lambda, \mu, \nu)$ is a critical point, which
leads to a contradiction as $\lambda, \mu, \nu$ are all positive.  If
$\eta$ is finite, then we also obtain a contradiction because
$(-\infty,\eta)$ is by definition a maximal interval on which the
solution exists.  So we deduce that $a$, and hence $\rho$, tends to
infinity as $u$ tends to $\eta$.

Therefore we must study the asymptotics of the metric $(\ref{metrho})$
as $\rho$ tends
 to infinity. It follows from $(\ref{eqrho})$ that
\[
\frac{d}{d\rho} \left( \rho^{2n+2}W \right)
 = 2n \rho^{2n+1} \left( \frac{a}{b}+
\frac{b}{a} \right) - \Lambda \rho^{2n+3}.
\]
Solving for $W$, and recalling that $\frac{b}{a}$ decreases monotonically
on $(-\infty, \eta)$ to some finite positive limit $L$, we see that
$W=O(\rho^{2})$ if $\Lambda$ is negative
and $W$ is bounded for large $\rho$ if $\Lambda$ is zero.
It follows that the geodesic distance
$\int_{\rho_{0}}^{\infty} \sqrt{W^{-1}} d \rho$
from $\rho=\rho_{0}$ to $\rho=\infty$ is infinite.

We have shown that the metric is complete. The underlying
topological manifold is the total space of
a rank $n+1$ vector bundle $E$ over $S^{n+1}$. In fact
the sphere bundle of this vector bundle is the Stiefel manifold
$SO(n+2)/SO(n)$, so $E$ is in fact the tangent bundle of
$S^{n+1}$. The Bolt is the zero section of $E$.

We summarise our results in the next theorem.
\begin{theorem} \label{theorem}
 The unstable curves of points $(0,K,K)$, where $K$ is nonzero, give complete
 Einstein metrics with nonpositive Einstein constant on $TS^{n+1}$.
\end{theorem}

\bigskip
\begin{remark}
(i) As remarked earlier, our metrics and the equations
 (\ref{KE1})-(\ref{KE3}) are invariant under
changing the sign of any of $a,b,c$, and also under interchanging $a,
b$.  Therefore the metrics arising from unstable curves of $(0,K,-K)$
or $(K, 0, \pm K)$ where $K$ is nonzero are isometric to those of
Theorem \ref{theorem}.

(ii) One can also obtain complete metrics by considering the unstable
curves of $(K, \pm K,0)$ where $K$ is nonzero and $n-\Lambda K^{2}$ is
a half-integer. The latter condition is needed to ensure that the
metric can be completed by adding a Bolt, which in this
case is the flag manifold $SO(n+2)/(SO(n) \times SO(2))$.
The underlying manifold of
the complete metric is the total space of a complex line bundle over
the Bolt.

 However, for these trajectories $a^{2}$ is identically equal to
$b^{2}$ and the metrics on the $SO(n+2)$ orbits are
obtained by Riemannian submersions with circle fibres.
The resulting cohomogeneity one metrics  are included in the examples of
B$\rm \acute{e}$rard-Bergery \cite{BB} (see also Gibbons and Pope
\cite{GP}, and Pedersen and Poon
\cite{P}, \cite{PP}).

(iii) If $\Lambda=0$ we have a special solution, given up to
translation of $u$ by
\[
a = b = \sqrt{\frac{n+1}{2n}} \; c =  \left(1 -
\left(\frac{2n}{n+1} \right)u \right)^{-1/2}
\]
This is the Ricci-flat cone over an Einstein metric on $SO(n+2)/SO(n)$.
For this solution, the trajectory is an unstable curve emanating from
$(0,0,0)$, hence the metric can be viewed as a limiting case of both
the Stenzel metrics and of the B\'erard-Bergery metrics, in the limit
that the volume of the Bolt ($S^{n+1}$ or $SO(n+2)/SO(n)
\times SO(2)$ respectively) tends to zero.

If $\Lambda < 0$, then, as remarked above, the values of $K$ for
which the B\'erard-Bergery metric extends over the Bolt are discrete
so we cannot take a continuous limit of the B\'erard-Bergery metrics.

(iv) The Bolt is complex for the B\'erard-Bergery metrics, but totally
real for the metrics arising from unstable curves of $(K, 0, K)$. We
can see this from (\ref{cplex}), as in the B\'erard-Bergery case the
tangent space to the Bolt is the $J$-invariant space $\p_1 \oplus
\p_2$, while in the other case the tangent space is $\p_1 \oplus
\p_3$, which satisfies $(\p_1 \oplus \p_3) \cap J(\p_1 \oplus \p_3) =
\{0 \}$.

(v) In the special case when $n=1$ and $\Lambda=0$ the equations
are symmetric in $a,b,c$ and the metrics in the two families are both
isometric to the hyperk\"ahler Eguchi-Hanson metric--but with respect to two
different complex structures in its two-sphere of complex structures.

\end{remark}

\section{Relation with Stenzel's work}

Stenzel \cite{St} has shown the existence
 of complete Ricci-flat K$\rm \ddot{a}$hler metrics
on the tangent bundles of spheres by different methods. He considers the
complex structure on $TS^{n+1}$ obtained by identifying this space
with the quadric in $\cl^{n+2}$. Regarding the complement of the zero section
in $TS^{n+1}$ as $(0, \infty) \times SO(n+2)/SO(n)$, and taking $r$ as a
coordinate on $(0, \infty)$, this complex structure is defined by
\[
\tilde{J} : \frac{\partial}{\partial r} \mapsto -X_{n+1,n+2}, \;\;\;\;
 X_{i,n+1} \mapsto -\tanh r \; X_{i,n+2}.
\]
Stenzel now shows that if $f$ satisfies the differential equation
\begin{equation}
\frac{d}{dr} \left( \left(\frac{df}{dr} \right)^{n+1} \right) =
(n+1)k (\sinh r)^{n} \label{eqstenzel1}
\end{equation}
for some constant $k$, and if the derivative of $f$ vanishes at $r=0$,
then the function $\phi$ defined by $\phi(r)=f(2r)$ is the K\"ahler
 potential with respect to $\tilde{J}$ for a Ricci-flat
metric $\tilde{g}$ on $(0, \infty) \times SO(n+2)/SO(n)$
which extends to the whole of $TS^{n+1}$.

Explicitly, the metric is
\begin{equation}
\tilde{g}=\frac{1}{2} \phi^{\prime \prime} d r^{2} + \tilde{g}_{r},
\label{metstenzel1}
\end{equation}
 where the homogeneous metric $\tilde{g}_{r}$
is defined in the notation of section 2 by
\begin{equation}
\frac{1}{2} \left( \phi^{\prime} (r) \tanh r \;B \mid_{\p_{1}} \oplus
\; \phi^{\prime} (r) \coth r \;B \mid_{\p_{2}} \oplus
\; \phi^{\prime \prime} (r) \;B \mid_{\p_{3}} \right). \label{metstenzel2}
\end{equation}
(In these expressions, and for the remainder of the paper, we use a prime
to indicate differentiation with respect to $r$).

By virtue of (\ref{eqstenzel1}) the function $\phi$ satisfies
\begin{equation}
\phi^{\prime \prime \prime} = n \phi^{\prime \prime} \left( \coth r +
\tanh r -\frac{\phi^{\prime \prime}}{\phi^{\prime}}
\right). \label{eqstenzel2}
\end{equation}

Defining the variable $u$ by $du= (2/\phi^{\prime}) dr$, it is now
straightforward to verify that this metric is one of those defined by
equations (\ref{KE1})-(\ref{KE3}) with $\Lambda=0$.  Moreover
$\phi^{\prime} \tanh r$ tends to zero as $r$ tends to zero. In order
to obtain a metric on $TS^{n+1}$ we need $\phi^{\prime} \coth r$ and
$\phi^{\prime \prime}$ to tend to equal nonzero limits as $r$ tends to
zero.  This implies that the solution to (\ref{KE1})-(\ref{KE3}) must
be an unstable curve of a critical point $(0,K, \pm K)$ with $K$
nonzero, so the metric is one of the Ricci-flat examples of
Theorem \ref{theorem}

Conversely, suppose we have a solution to the equations (\ref{KE1})-(\ref{KE3})
 with $\Lambda=0$ arising as an unstable curve of $(0,K,\pm K)$ where
 $K$ is nonzero. Letting $r$ satisfy $dr=ab \;du$, we obtain from
 (\ref{KE1}),(\ref{KE2}) the equations
\[
(ab)^{\prime} =c^{2}, \;\;\; : \;\;\;
\left(\frac{a}{b} \right)^{\prime} = 1 - \left(\frac{a}{b} \right)^{2}.
\]
It easily follows that the metric can be put in the form given by
(\ref{metstenzel1}),(\ref{metstenzel2}).  Equation (\ref{KE3}) is then
equivalent to (\ref{eqstenzel2}), which may be integrated to give
(\ref{eqstenzel1}). Moreover, the fact that $(a,b,c)$ tends to $(0,K,
\pm K)$ as $u$ tends to $-\infty$ gives the initial condition for
equation (\ref{eqstenzel1}).  Therefore the Ricci-flat metrics of
Theorem \ref{theorem} are precisely the metrics of \cite{St}

This argument also shows that any of the metrics of Theorem
\ref{theorem} with negative Einstein constant can be put in the
form (\ref{metstenzel1}), (\ref{metstenzel2}). Now the coordinate
$r$ is defined by $r = \int_{-\infty}^{u} ab \;\;du$ or
equivalently, in the notation of (\ref{eqrho}), $r=
\int_{0}^{\rho(u)} \frac{2}{\rho W} d \rho$.  We saw earlier that
$W=O(\rho^{2})$ if the Einstein constant is negative, so this
integral converges as $\rho$ tends to $\infty$, that is, as $u$
tends to $\eta$. The upshot is that the complex structure for our
complete metrics with $\Lambda$ negative is obtained by
restricting Stenzel's complex structure on $T S^{n+1}$($=$ quadric
in $\cl^{n+2}$) to a tube of finite radius.  In particular the
K$\rm \ddot{a}$hler structure extends to the Bolt, giving us a
globally defined K$\rm \ddot{a}$hler structure in the metrics of
Theorem \ref{theorem}.

\bigskip
{\bf Acknowledgements} A.S.D. thanks the Max-Planck-Institut f\"ur
Mathematik, Bonn, for its support, and Xenia de la Ossa and Robert
Sz\"oke for discussions. I.A.B.S. would like to thank the University
of Newcastle for a Wilfred Hall fellowship.

\begin{flushleft}
{\small Dancer : Jesus College, Oxford University, Oxford, OX1 3DW,
United Kingdom} \\{\tt Email : dancer@maths.ox.ac.uk}

\medskip
{\small Strachan : Department of Mathematics, University of Hull,
Hull, HU6 7RX, United Kingdom} \\{\tt Email :
i.a.strachan@hull.ac.uk}
\end{flushleft}

\end{document}